\documentclass{amsart}
\usepackage{amssymb}

\newtheorem{theorem}{Theorem}[section]
\newtheorem{conjecture}[theorem]{Conjecture}
\newtheorem{corollary}[theorem]{Corollary}
\newtheorem{lemma}[theorem]{Lemma}
\newtheorem{proposition}[theorem]{Proposition}
\newtheorem{definition}[theorem]{Definition}
\begin{document}
\title{Note on the conjecture of D.~Blair in contact Riemannian geometry.}
\author{Vladimir Krouglov.}
\begin{abstract}
The conjecture of D.~Blair says that there are no nonflat Riemannian metrics of nonpositive curvature compatible with a contact structure. We prove this conjecture for a certain class of contact structures on closed $3$-dimensional manifolds and construct a local counterexample. We also prove that a hyperbolic metric on $\mathbb{R}^3$ cannot be compatible with any contact structure. \end{abstract}
\maketitle

\section{Introduction}
\noindent In \cite[p.~99]{DB} author states the following conjecture:
\begin{conjecture}
There are no nonflat Riemannian metrics of nonpositive curvature that are compatible with a contact structure. 
\end{conjecture}
On a $3$-torus, standard Euclidean metric is compatible with a contact structure given by the kernel of the one-form $\cos(z)dx + \sin(z)dy$. Despite the fact that there exist contact structures on higher dimensional tori \cite{B}, as it was shown in \cite{DB2} flat metric cannot be compatible with any contact structure when the dimension of a manifold is greater than $3$. 

Using the result of A.~Zeghib \cite{Z} on the existence of geodesic flows, on closed manifolds the conjecture of D.~Blair is true for the Riemannian metrics of strictly negative sectional curvature. In \cite{L}, it has been shown that the conjecture is true for the homogenous Riemannian metric adapted to a homogenous contact structure. 

Note, that in view of the results in \cite{EKM} closed contact metric manifolds of nonpositive curvature would provide a source of examples of tight contact structures.

The main result of the present paper is the proof of the D.~Blair's conjecture for contact structures which are sufficiently nontrivial as fibrations. We prove the following
\begin{theorem}
Assume that $M$ is a closed $3$-manifold with a contact structure $\xi$ which cannot be decomposed as a sum of two one-dimensional fibrations
$\xi \ne \eta_1 \oplus \eta_2$. Then the conjecture of D.~Blair is true for $(M, \xi)$.
\end{theorem}
By the result of Z.~Olszak in \cite{O}, when the dimension of a manifold is greater than three constant negative curvature metrics cannot be compatible with a contact structure (even when the manifold is not compact). Analyzing the curvature tensor of a compatible metric we prove this result in dimension three. 
\begin{proposition}
Constant negative curvature metric on a $3$-manifold cannot be compatible with a contact structure.
\end{proposition}
We end with a local counterexample to the conjecture of D.~Blair. We construct a Riemannian metric compatible with a standard contact structure on $\mathbb{R}^3$ which has strictly negative curvature in some neighborhood of zero in $\mathbb{R}^3$.
\section{Contact metric manifolds.}
\subsection{Compatible metrics.}
Assume that $(M, \xi)$ is a contact $3$-manifold. If we fix a one-form $\alpha$ among the conformal class $\{f\alpha' : \  \mbox{for positive functions $f$ on $M$}\}$ which we call the contact one-form associated with the contact structure then there is a unique vector field $N$ called the Reeb vector field of $\alpha$ such that
$$
\alpha(N)=1, \ L_N\alpha=\iota_Nd\alpha=0
$$
Let $J$ be an almost complex structure on $\xi$ (i.e. $J^2=-id$). We may complement it to a linear operator on $TM$ by setting $JN=0$.
\begin{definition}
A Riemannian metric $\langle \cdot, \cdot \rangle$ is called compatible with $\xi$ if there is an associated $1$-form $\alpha$ and an almost complex structure $J$ such that
$$
\langle N, X \rangle = \alpha(X), \ k \langle X, JY \rangle = d\alpha(X, Y)
$$
where $k$ is some constant and $X$ and $Y$ are the vector fields on $M$.

By a contact metric manifold we are going to understand the tuple $(M, \xi, \alpha, \langle \cdot, \cdot \rangle, J)$.
\end{definition}
\subsection{Second fundamental form.}
The second fundamental form of a plane field is a symmetric bilinear form which generalizes the corresponding notion for a surface inside the Riemannian manifold. The following definition is due to Reinhart \cite{Re}
\begin{definition}
The second fundamental form of plane field $\xi$ is a bilinear form on $\xi$ defined as
$$
II(X, Y) = \frac{1}{2}\langle\nabla_X Y + \nabla_Y X, N\rangle
$$
where $X$ and $Y$ are in $\xi$, $N$ is a unit normal vector field to $\xi$ and $\nabla$ is a Levi-Civita connection of $\langle \cdot, \cdot \rangle$.
\end{definition}
We are going to call the linear operator $A_N$ which corresponds to $II$ with respect to $\langle \cdot, \cdot \rangle$ -- a \emph{shape operator} of $\xi$. Since $II$ is symmetric, the shape operator has two real eigenvalues that we call the principal curvatures of $\xi$. The eigenvectors of $A_N$ will be called the principal directions of $\xi$. We also define the \emph{extrinsic curvature} $K_e$ and the \emph{mean curvature} $H$ of $\xi$ as the determinant and the half trace of the \emph{shape operator} correspondingly. When the plane field $\xi$ is integrable, the second fundamental form of $\xi$ coincides with a second fundamental forms of the integral surfaces. All notions of the classic surface theory extend naturally to the context of plane distributions.
\subsection{Extrinsic geometry in compatible metric.}
When $M$ is a contact metric manifold, the contact structure $\xi$ has a very special geometry with respect to the compatible metric $\langle \cdot, \cdot \rangle$. We have the following
\begin{proposition}\cite{DB}
With respect to a compatible metric, the Reeb vector field is a unit speed geodesic vector field and the contact structure is minimal.
\end{proposition}
We are also going to summarize several properties of the contact structures with respect to a compatible metric that will be used in the derivation of the curvature tensor.
\begin{lemma}\label{TL}
Let $(M, \xi, \alpha,\langle \cdot, \cdot \rangle, J)$ be a contact metric manifold. Then,
\begin{enumerate}
\item{$J$ is a rotation by $\frac{\pi}{2}$ in $\xi$.}
\item{For every pair of orthonormal vectors $X$ and $Y$ in $\xi$ the function $\langle [X, Y], N \rangle=\pm k$.}
\item{If $X$ and $Y$ are unit orthogonal principal directions of $A_N$ then $$\langle \nabla_X Y, N \rangle = - \langle \nabla_Y X, N \rangle = \pm \frac{k}{2}.$$ }
\end{enumerate}
\end{lemma}
\noindent \emph{Proof:} For every pair of vectors $X$ and $Y$ in $\xi$
$$
k\langle JX, JY \rangle = d\alpha(JX, Y) = - d\alpha(Y, JX) = - k\langle Y, J^2Y \rangle = k\langle X, Y \rangle.
$$
We are left to check that $X$ is orthogonal to $JX$. This follows from
$$
k\langle X, JX \rangle = d\alpha(X, X) = 0
$$
If $X$ and $Y$ are orthonormal, then $Y = \pm JX$. We have
$$
d\alpha(X, Y) = X\alpha(Y) - Y\alpha(X)-\alpha([X, Y])=-\langle[X, Y], N \rangle.
$$
On the other hand
$$
d\alpha(X, Y) = \pm k\langle X, X \rangle = \pm k
$$
which proves $(2)$.

Since $X$ and $Y$ are the eigenvectors of $A_N$, $\frac{1}{2}\langle \nabla_X Y + \nabla_Y X, N \rangle=0$. From $(2)$,
$$
\langle \nabla_X Y, N \rangle = \frac{1}{2}\langle \nabla_X Y + \nabla_Y X, N \rangle + \frac{1}{2}\langle \nabla_X Y - \nabla_Y X, N \rangle = \pm \frac{k}{2}
$$

\section{Curvature tensor of the compatible metric on a $3$-manifold.}
In this section we are going to compute the matrix of the curvature tensor of a compatible metric. Assume that $(M, \xi, \alpha, \langle \cdot, \cdot \rangle, J)$ is a contact metric manifold. Let $N$ be the Reeb vector field of $\alpha$. Denote by $X$ and $Y$ the (local) orthonormal frame in $\xi$ that consists of the eigenvectors of the shape operator at a given point $p \in M$.

Let $\lambda$ be a principal curvature that corresponds to a principal direction $X$. Since $\xi$ is minimal, the mean curvature of $\xi$ vanishes and $Y$ corresponds to the principal curvature $-\lambda$.

\begin{lemma}
With respect to a basis of bivectors $X \wedge Y$, $X \wedge N$ and $Y \wedge N$ the matrix of the curvature tensor of $\langle \cdot, \cdot \rangle$ is given by
$$\mathcal{R}=\left(
\begin{array}{ccc}
-\frac{3k^2}{4} + \lambda^2 + K  & -Y(\lambda) - 2\lambda \langle \nabla_X X, Y \rangle & X(\lambda) - 2\lambda \langle \nabla_Y Y, X \rangle\\
-Y(\lambda) - 2\lambda \langle \nabla_X X, Y \rangle & \frac{k^2}{4} - \lambda^2 + N(\lambda) &  2\lambda \langle \nabla_N X,  Y \rangle\\
X(\lambda) - 2\lambda \langle \nabla_Y Y, X \rangle &  2\lambda \langle \nabla_N X,  Y  \rangle & \frac{k^2}{4} - \lambda^2 + N(\lambda)
\end{array}
\right)$$
where
$$
K = X (\langle \nabla_Y Y, X \rangle) + Y (\langle \nabla_X X, Y \rangle) - \langle \nabla_Y Y, X \rangle^2 - \langle \nabla_X X, Y \rangle^2 - \langle [X, Y], N \rangle \langle [N, Y], X \rangle
$$
is the curvature of a generalized Webster connection (see \cite{T} for the definition) and $\lambda$ is an eigenvalue of the shape operator which corresponds to $X$.
\end{lemma}
\noindent \emph{Proof:} By replacing $X$ by $-X$ if required we may assume that $\langle [X, Y], N \rangle=k$. \\
\textbf{Calculation of $\mathcal{R}_{11} = \langle R(X, Y)Y, X\rangle$}.
$$
\langle R(X, Y) Y, X\rangle = \langle \nabla_X \nabla_Y Y, X \rangle - \langle \nabla_Y \nabla_X Y, X \rangle - \langle \nabla_{[X, Y]} Y, X \rangle
$$
The first summand is
\begin{align*}
\begin{split}
\langle \nabla_X \nabla_Y Y, X \rangle &= X (\langle \nabla_Y Y, X \rangle) - \langle \nabla_Y Y,  \nabla_X X \rangle \\
&=X (\langle \nabla_Y Y, X \rangle) - \langle \nabla_Y Y, N \rangle \langle \nabla_X X, N \rangle = X (\langle \nabla_Y Y, X \rangle) + \lambda^2
\end{split}
\end{align*}
The second summand is
$$
- \langle \nabla_Y \nabla_X Y, X \rangle = - Y (\langle \nabla_X Y, X \rangle) + \langle  \nabla_X Y, \nabla_Y X \rangle  = Y (\langle \nabla_X X, Y \rangle) - \frac{k^2}{4}
$$
as follows from $(3)$ in Lemma \ref{TL}. The third summand is
\begin{align*}
\begin{split}
- \langle \nabla_{[X, Y]} Y, X \rangle &= -\langle [X, Y], X \rangle \langle \nabla_X Y, X \rangle - \langle [X, Y], Y \rangle \langle \nabla_Y Y, X \rangle - \langle [X, Y], N \rangle \langle \nabla_N Y, X \rangle \\
&= - \langle \nabla_X X, Y \rangle^2 - \langle \nabla_Y Y, X \rangle^2 - \langle [X, Y], N \rangle (\langle \nabla_Y N, X \rangle + \langle [N, Y], X \rangle)\\
&= - \langle \nabla_X X, Y \rangle^2 - \langle \nabla_Y Y, X \rangle^2 - \frac{k^2}{2} - \langle [X, Y], N \rangle \langle [N, Y], X \rangle
\end{split}
\end{align*}
Summing this up will give us the desired expression for $\mathcal{R}_11$.  \\
\textbf{Calculation of $\mathcal{R}_{22} = \langle R(X, N)N, X\rangle$}.
$$
\langle R(N, X) X, N\rangle = \langle \nabla_N \nabla_X X, N \rangle - \langle \nabla_X \nabla_N X, N \rangle - \langle \nabla_{[N, X]} X, N \rangle
$$
The first summand is
$$
\langle \nabla_N \nabla_X X, N \rangle =  \nabla_N \langle \nabla_X X, N \rangle - \langle \nabla_X X,  \nabla_N N \rangle = N(\lambda)
$$
The second summand is
\begin{align*}
\begin{split}
- \langle \nabla_X \nabla_N X, N \rangle &= -X \langle \nabla_N X, N \rangle + \langle \nabla_N X,  \nabla_X N \rangle \\
&= -X (N \langle X, N \rangle - \langle X, \nabla_N N \rangle) + \langle \nabla_N X, \nabla_X N \rangle = \langle \nabla_N X, \nabla_X N \rangle.
\end{split}
\end{align*}
Here we used that $N$ is a geodesic vector field. Finally, the last summand is
\begin{align*}
\begin{split}
- \langle \nabla_{[N, X]} X, N \rangle &= - \langle [N, X], X \rangle \langle \nabla_X X, N\rangle + - \langle [N, X], Y \rangle \langle \nabla_Y X, N\rangle \\
 &= - \lambda^2 - \langle [N, X], Y \rangle \langle \nabla_Y X, N\rangle
\end{split}
\end{align*}
Summing these expressions we get
$$
K(N ,X) = N(\lambda) - \lambda^2 - \langle [N, X], Y \rangle \langle \nabla_Y X, N\rangle + \langle \nabla_N X, Y \rangle \langle Y, \nabla_X N \rangle
$$
Using $(2)$ and $(3)$ of Lemma \ref{TL} we get
$$
K(N ,X) = N(\lambda) - \lambda^2 - \frac{k}{2}(\langle [N, X], Y \rangle  - \langle \nabla_N X, Y \rangle) = N(\lambda) - \lambda^2 +\frac{k^2}{4}
$$
\textbf{Calculation of $\mathcal{R}_{33} = \langle R(Y, N)N, Y\rangle$}. \\
By exactly the same calculations replacing $X$ by $Y$ we get
$$
\langle R(Y, N)N, Y\rangle = -N(\lambda) - \lambda^2 + \frac{k^2}{4}.
$$
\textbf{Calculation of $\mathcal{R}_{23}=\langle R(X, N) N, Y \rangle$}.
$$
\langle R(X, N) N, Y \rangle = \langle \nabla_X \nabla_N N, Y \rangle - \langle \nabla_N \nabla_X N, Y \rangle - \langle \nabla_{[X, N]} N, Y \rangle
$$
Obviously, since $N$ is geodesic the first summand is zero. Rewrite the second summand,
\begin{align*}
\begin{split}
- \langle \nabla_N \nabla_X N, Y \rangle &= - N(\langle \nabla_X N, Y \rangle) + \langle \nabla_X N,  \nabla_N Y \rangle = \langle \nabla_X N,  \nabla_N Y \rangle \\
&= \langle \nabla_X N, X \rangle \langle X, \nabla_N Y \rangle + \langle \nabla_X N, Y \rangle \langle Y, \nabla_N Y \rangle = - \lambda \langle X, \nabla_N Y \rangle
\end{split}
\end{align*}
The last summand
$$
-\langle \nabla_{[X, N]} N, Y \rangle = -\langle [X, N], X \rangle \langle \nabla_X N, Y \rangle  - \langle [X, N], Y \rangle \langle \nabla_Y N, Y \rangle
$$
Summing these expressions we get:
$$
\langle R(Y, N)N, Y\rangle = \lambda \langle \nabla_N X,  Y \rangle + \lambda \langle \nabla_X N, Y \rangle  - \lambda \langle [X, N], Y \rangle = 2\lambda \langle \nabla_N X,  Y \rangle
$$
\textbf{Calculation of $\mathcal{R}_{13}=\langle R(X, Y) N, Y \rangle$}.
$$
\langle R(X, Y) N, Y \rangle = -\langle R(X, Y) Y, N \rangle = -\langle \nabla_X \nabla_Y Y, N \rangle + \langle \nabla_Y \nabla_X Y, N \rangle + \langle \nabla_{[X,Y]} Y, N \rangle
$$
The first summand is
\begin{align*}
\begin{split}
-\langle \nabla_X \nabla_Y Y, N \rangle &= -X \langle \nabla_Y Y, N \rangle + \langle  \nabla_Y Y, \nabla_X N \rangle =  X(\lambda) + \langle  \nabla_Y Y, X \rangle \langle X, \nabla_X N \rangle \\
&= X(\lambda) - \lambda \langle \nabla_Y Y, X \rangle
\end{split}
\end{align*}
The second summand is
\begin{align*}
\begin{split}
\langle \nabla_Y \nabla_X Y, N \rangle &= Y(\langle \nabla_X Y, N \rangle) - \langle  \nabla_X Y, \nabla_Y N \rangle = -\langle  \nabla_X Y, X \rangle \langle X, \nabla_Y N \rangle\\ &= -\langle  \nabla_X X, Y \rangle \langle N, \nabla_Y X \rangle
\end{split}
\end{align*}
Finally, the last summand,
\begin{align*}
\begin{split}
\langle \nabla_{[X,Y]} Y, N \rangle &= \langle [X, Y], X \rangle \langle \nabla_X Y, N \rangle + \langle [X, Y], Y \rangle \langle \nabla_Y Y, N \rangle\\ &= -\langle  \nabla_X X, Y \rangle \langle N, \nabla_X Y \rangle -\lambda \langle \nabla_Y Y, X \rangle
\end{split}
\end{align*}
Summing this up gives us
$$
\langle R(X, Y) N, Y \rangle=X(\lambda) - 2\lambda \langle \nabla_Y Y, X \rangle
$$
\textbf{Calculation of $\mathcal{R}_{12}=\langle R(X, Y) N, X \rangle$}. \\
Analogously,
$$
\langle R(X, Y) N, X \rangle=-Y(\lambda) - 2\lambda \langle \nabla_X X, Y \rangle
$$

\begin{corollary}
Assume that $M$ is a closed $3$-manifold with a contact structure $\xi$ which cannot be decomposed as a sum of two one-dimensional fibrations
$\xi \ne \eta_1 \oplus \eta_2$. Then the conjecture of D.~Blair is true for $(M, \xi)$.
\end{corollary}
\noindent \emph{Proof:} Under the assumptions of the corollary, for every Riemannian metric $g$ on $M$, $\xi$ must have an \emph{umbilic point}. At this point we have $\lambda=0$ and
$$
K(X, N) + K(Y, N) = \frac{k^2}{2} - 2\lambda^2 = \frac{k^2}{2} > 0
$$
Therefore, $g$ cannot have nonpositive curvature.
\begin{proposition}
Constant negative curvature metric cannot be compatible with any contact structure.
\end{proposition}
\noindent \emph{Proof:} Since the sectional curvature of the metric is constant, it is easy to see that the principal curvatures should also be constant on $M$. With respect to a basis of bivectors $X \wedge Y$, $X \wedge N$ and $Y \wedge N$ the matrix of the curvature operator of $\langle \cdot, \cdot \rangle$ is diagonal. Noting that $\lambda$ cannot be zero,
$$
\left\{\begin{array}{l} \langle \nabla_N X, Y \rangle = 0 \\
\langle \nabla_X X, Y \rangle = 0 \\
\langle \nabla_Y Y, X \rangle = 0
\end{array}
\right.
$$
The Webster curvature with respect to $\langle \cdot, \cdot \rangle$ is
\begin{align*}
\begin{split}
K =& X (\langle \nabla_Y Y, X \rangle) + Y (\langle \nabla_X X, Y \rangle) - \langle \nabla_Y Y, X \rangle^2 \\ &- \langle \nabla_X X, Y \rangle^2 - \langle [X, Y], N \rangle \langle [N, Y], X \rangle \\
=& -\frac{k}{2} (\langle \nabla_N Y, X \rangle - \langle \nabla_Y N, X \rangle) =\frac{k}{2} ( \langle \nabla_Y N, X \rangle) = \frac{k^2}{2}
\end{split}
\end{align*}
Since $K(X, Y) = K(X, N)$ we have that
$$
-\frac{3k^2}{4}+ \frac{k^2}{2} + \lambda^2 = \frac{k^2}{4}-\lambda^2,
$$
which implies that $\lambda = \frac{1}{2}$ and the metric is flat.
\section{Local counterexample to the conjecture of D.~Blair.}
On $\mathbb{R}^3$ with cartesian coordinates $(x,y,z)$ consider a standard contact structure $\xi$ given by the kernel of the one-form $\alpha=dz-xdy$. We will construct a Riemannian metric which would be compatible with $\xi$ and have nonpositive (even strictly negative) curvature in some neighborhood of zero in $\mathbb{R}^3$.

With respect to this metric the Reeb vector field of $\alpha$ has to be a unit geodesic vector field and $\xi$ has to be a minimal distribution. It is easy to check that in this case the matrix of $g$ should have the form
$$
g = \left(\begin{array}{ccc}a & b & 0 \\ b & c & x \\ 0 & x & 1 \end{array}\right).
$$
where the functions $a, b$ and $c$ additionally satisfy the condition
$$
H = \frac{1}{2}\frac{\partial}{\partial z}(a(c-x^2)-b^2) = 0
$$
This condition will be automatically satisfied if we choose
$$\left\{
\begin{array}{l}
a = Ae^z \\
b = 1 \\
c = x^2 + Be^{-z}
\end{array}
\right.$$
With respect to an orthonormal frame $(\frac{\partial}{\partial z}, \frac{1}{\sqrt{Ae^z}}\frac{\partial}{\partial x}, \sqrt{\frac{Ae^z}{AB-1}}(-\frac{1}{Ae^z}\frac{\partial}{\partial x} + \frac{\partial}{\partial y}-x\frac{\partial}{\partial z}))$ the curvature tensor is given by a matrix
$$
\left( \begin{array}{ccc}
\frac{1}{4}\frac{AB-3-2x^2Ae^z}{AB-1} & -\frac{1}{2}x\sqrt{\frac{Ae^z}{AB-1}} & \frac{1}{2}\frac{x\sqrt{Ae^z}}{AB-1} \\
-\frac{1}{2}x\sqrt{\frac{Ae^z}{AB-1}} & -\frac{1}{4} & 0 \\
\frac{1}{2}\frac{x\sqrt{Ae^z}}{AB-1} & 0 & \frac{1}{4}
\end{array}\right)
$$
Clearly when $AB \in (1,3)$, the matrix of the curvature tensor is negatively definite in some neighborhood of zero in $\mathbb{R}^3$.

\end{document}